\documentclass[preprint]{elsarticle}

\usepackage{amsmath}
\usepackage{amsthm}
\usepackage{amssymb}
\usepackage{bm}

\newcommand{\grad}{\mathop{\rm grad}\nolimits}
\renewcommand{\div}{\mathop{\rm div}\nolimits}

\begin{document}

\begin{frontmatter}

\title{Numerical solving the identification problem for the lower coefficient of parabolic equation\tnoteref{label1}}
\tnotetext[label1]{This work was supported by RFBR (project 13-01-00719)}

\author{P.N. Vabishchevich\corref{cor1}\fnref{lab1}}
\ead{vabishchevich@gmail.com}
\cortext[cor1]{Correspondibg author.}

\author{V.I. Vasil'ev\fnref{lab2}}
\ead{vasvasil@mail.ru}

\address[lab1]{Nuclear Safety Institute, Russian Academy of Sciences,
              52, B. Tulskaya, 115191 Moscow, Russia}

\address[lab2]{North-Eastern Federal University,
	      58, Belinskogo, 677000 Yakutsk, Russia}

\begin{abstract}
In the theory and practice of inverse problems for partial differential equations (PDEs) much attention is paid to the problem of the identification of coefficients from some additional information.
This work deals with the problem of determining in a multidimensional parabolic equation the lower coefficient that depends on time only. To solve numerically a nonlinear inverse problem, linearized approximations in time are constructed using standard finite element procedures in space.
The computational algorithm is based on a special decomposition, where the transition to a new time level is implemented via solving two standard elliptic problems. The numerical results presented here for a model 2D problem demonstrate  capabilities of the proposed computational algorithms for approximate solving inverse problems.
\end{abstract}


\end{frontmatter}

\section{Introduction} 
\label{sec:1}

Mathematical modeling of many applied problems of science and engineering results in the numerical solution of inverse problems \cite{alifanov2011inverse,tarantola1987inverse,aster2011parameter}.
Inverse problems often belong to the class of ill-posed (conditionally correct) problems, and therefore various regularization algorithms are employed to solve them numerically \cite{kirsch1996introduction,tikhonov1977solutions,engl2000regularization}. 

Particular attention should be given to inverse problems for PDEs \cite{lavrentev1986ill,isakov1998inverse}.
In this case,  a theoretical study includes the fundamental questions of uniqueness of the solution and its stability both from the viewpoint of the theory of differential equations
\cite{express2000methods,belov2002inverse} and from the viewpoint of  the theory of optimal control  for distributed systems \cite{maksimov2002dynamical}.
Many inverse problems are formulated as non-classical problems for PDEs.
To solve these problems approximately, emphasis is on  the development of stable computational algorithms that take into account peculiarities of  inverse problems \cite{vogel2002computational,samarskii2007numerical}.
 
Among inverse problems for PDEs we distinguish coefficient inverse problems, which are associated with the identification of coefficients and/or the right-hand side of an equation using some additional information. When considering  time-dependent problems, the identification of the coefficient dependences on space and on time is usually separated into individual problems \cite{isakov1998inverse,express2000methods}.
In some cases, we have  linear inverse problems (e.g., identification problems for the right-hand side of an equation); this situation essentially simplify their study.

Much attention is paid to the problem of determining the lower coefficient of a parabolic equation of second order, where, in particular, the coefficient depends on time only.
An additional condition is most often formulated as a specification of the solution at an interior point or  as the average value integrated over the whole domain.
The existence and uniqueness of the solution of such an inverse problem and well-posedness of this problem in various functional classes are examined, for example, in the works \cite{prilepko1985determination,prilepko1987solvability,macbain1986existence,cannon1990inverse}.

Numerical methods for solving the problem of the identification of the lower coefficient of parabolic equations are considered in many works \cite{wang1989finite,cannon1994numerical,dehghan2005identification,ye2007stability,wang2012inverse}.
In view of the practical use, we highlight separately studies dealing with numerical solving  inverse problems for multidimensional parabolic equations \cite{dehghan2001numerical,daoud2005splitting}.
To construct computational algorithms for the identification of the lower coefficient of a parabolic equation, there is widely used  the idea of transformation of the equation by introducing new unknowns that results in a linear inverse problem.
 
In this paper,  for a multidimensional parabolic equation, we consider the problem of determining the lower coefficient that depends on time only.  Approximation in space is performed using standard finite elements
\cite{quarteroni2008numerical,brenner2008mathematical}.
The main features of the nonlinear inverse problem are taken into account via a proper choice of the linearized approximation in time.
Linear problems at a particular time level are solved on the basis of a special decomposition into two standard elliptic problems. The paper is organized as follows. 
In Section~\ref{sec:2}, for a parabolic equation of second order, we formulate the inverse problem of the identification of the lower coefficient.
The computational algorithm  based on the linearization scheme  is described in Section~\ref{sec:3}.
Section~\ref{sec:4} presents possibilities of the schemes with the second-order approximation in time. Numerical results for a model 2D  inverse problem are  discussed in Section~\ref{sec:5}.

\section{Problem formulation}
\label{sec:2}

For simplicity, we restrict ourselves to a 2D problem. Generalization to the 3D case is trivial. 
Let ${\bm x} = (x_1, x_2)$ and $\Omega$ be a  bounded polygon.
The direct problem is formulated as follows.
We search $u({\bm x},t)$, $0 \leq t \leq T, \ T > 0$ such that 
it is the solution of the parabolic equation of second order:
\begin{equation}\label{eq:1}
  \frac{\partial u}{\partial t}- \div (k({\bm x}) \grad u) + p(t) u = f({\bm x},t),
  \quad {\bm x} \in \Omega,
  \quad 0 < t \leq T .    
\end{equation} 
The boundary and initial conditions are also specified:
\begin{equation}\label{eq:2}
  k({\bm x}) \frac{\partial u}{\partial n} + g({\bm x}) u = 0,
  \quad {\bm x} \in \partial\Omega,
  \quad 0 < t \leq T,    
\end{equation} 
\begin{equation}\label{eq:3}
  u({\bm x}, 0) = u_0({\bm x}),
  \quad {\bm x} \in \Omega ,  
\end{equation} 
where $n$ is the normal to $\Omega$.
The formulation (\ref{eq:1})--(\ref{eq:3}) presents the direct problem, where the right-hand side, coefficients of the equation as well as  the boundary and initial conditions are specified.

Let us consider the inverse problem, where in equation (\ref{eq:1}), the coefficient
$p(t)$ is unknown. An additional condition is often formulated as 
\begin{equation}\label{eq:4}
  \int_{\Omega} u({\bm x}, t) \omega({\bm x}) d {\bm x}  = \varphi(t),
  \quad 0 < t \leq T,    
\end{equation} 
where $\omega({\bm x})$ is a weight function.
In particular, choosing $\omega({\bm x}) = \delta({\bm x} -{\bm x}^*)$ (${\bm x}^* \in \Omega$),
where $\delta({\bm x})$ is the Dirac $\delta$-function, from (\ref{eq:4}), we get
\begin{equation}\label{eq:5}
  u({\bm x}^*, t) = \varphi(t),
  \quad 0 < t \leq T .   
\end{equation} 
 
We assume that the above inverse problem of finding a pair of $u({\bm x},t), \ p(t)$ from equations  
(\ref{eq:1})--(\ref{eq:3}) and additional conditions (\ref{eq:4}) or
(\ref{eq:5}) is well-posed. The corresponding conditions for existence and uniqueness of the solution are available in the above-mentioned works.
In this paper, we  consider only the numerical solution of  these inverse problems omitting theoretical issues of the convergence of an approximate solution to the exact one.

From the nonlinear inverse problem we can proceed to the linear one.
Suppose
\[
  v({\bm x}, t) = \chi(t) u({\bm x}, t),
  \quad  \chi(t) = \exp \left (\int_{0}^{t} p(\theta) d \theta \right) .  
\] 
Then from (\ref{eq:1})--(\ref{eq:3}), we get
\[
  \frac{\partial v}{\partial t}- \div (k({\bm x}) \grad v)  = \chi(t) f({\bm x},t),
  \quad {\bm x} \in \Omega,
  \quad 0 < t \leq T ,  
\] 
\[
  k({\bm x}) \frac{\partial v}{\partial n} + g({\bm x}) v = 0,
  \quad {\bm x} \in \partial\Omega,
  \quad 0 < t \leq T,    
\]
\[
  v({\bm x}, 0) = u_0({\bm x}),
  \quad {\bm x} \in \Omega .  
\]
The additional conditions (\ref{eq:4}) and (\ref{eq:5}) to  identify uniquely
$v({\bm x},t), \ \chi(t)$  take the form
\[
  \int_{\Omega} v({\bm x}, t) \omega({\bm x}) d {\bm x}  = \chi(t)\varphi(t),
  \quad 0 < t \leq T,   
\] 
\[
  u({\bm x}^*, t) =  \chi(t) \varphi(t),
  \quad 0 < t \leq T .   
\]
The above transition from the nonlinear inverse problem to the linear one
is in common use for numerical solving problems of identification.
In our work, we focus on the original formulation of the inverse problem
(\ref{eq:1})--(\ref{eq:4}) (or (\ref{eq:1})--(\ref{eq:3}), (\ref{eq:5}))
without going to the linear problem.

\section{Computational algorithm} 
\label{sec:3}

The inverse problem of determining the pair of $u({\bm x},t), \ p(t)$ is nonlinear.
The standard approach is based on the simplest approximations in time and involves the
iterative solution of the corresponding nonlinear problem for the evaluation
of the approximate solution at a new level. In our work, we apply such
approximations in time  that lead to linear problems for
evaluating the solution at the new time level.

Let us define a uniform grid in time  
\[
  \overline{\omega}_\tau =
  \omega_\tau\cup \{T\} =
  \{t^n=n\tau,
  \quad n=0,1,...,N,
  \quad \tau N=T\} 
\]
and denote $y^n = y(t^n), \ t^n = n \tau$.
Finite element approximations in space are employed.
In the polygon $\Omega$, we perform a triangulation and introduce for this computational grid 
a finite-dimensional space $V \subset H^1(\Omega)$ of finite elements.

Using the fully implicit scheme for approximation in time, we obtain
the following variational problem:
\begin{equation}\label{eq:6}
\begin{split}
  \int_{\Omega} \frac{u^{n+1} - u^n}{\tau} v  d {\bm x} & +
  \int_{\Omega} k({\bm x}) \grad u^{n+1} \grad v d {\bm x} \\
  & +
  \int_{\partial \Omega} g({\bm x}) u^{n+1} v dx +
  p^{n+1} \int_{\Omega} u^{n+1} v d {\bm x} \\
  & = 
  \int_{\Omega} f({\bm x},t^{n+1}) v d {\bm x},
  \quad \forall v \in V,   
  \quad n = 0,1, ..., N-1, 
\end{split}
\end{equation} 
\begin{equation}\label{eq:7}
  \int_{\Omega} u^0 v d {\bm x} = 
  \int_{\Omega} u_0 v d {\bm x} .
\end{equation}
The additional relations (\ref{eq:4}) and (\ref{eq:5}) take the form
\begin{equation}\label{eq:8}
  \int_{\Omega} u^{n+1} \omega({\bm x})  d {\bm x} = \varphi^{n+1}, 
\end{equation} 
\begin{equation}\label{eq:9}
  u^{n+1}({\bm x}^*) = \varphi^{n+1},
  \quad n = 0,1, ..., N-1 .  
\end{equation} 

To evaluate at the new time level the approximate solution $u^{n+1}({\bm x}), \ p^{n+1}$ from
(\ref{eq:6})--(\ref{eq:8}) or (\ref{eq:6}), (\ref{eq:7}), (\ref{eq:9}), 
some iterative procedures are necessary.
In solving  time-dependent problems, the solution slightly varies when it  pass
from the previous time level to the next one.
This basic feature of time-dependent problems is widely used in numerical
solving nonlinear problems through the application of linearization procedures.
We use a similar approach for the numerical solution of the
inverse problem that is concerned with the identification of the lower coefficient of a parabolic equation.

Instead of (\ref{eq:6}), we will solve the following equation:
\begin{equation}\label{eq:10}
\begin{split}
  \int_{\Omega} \frac{u^{n+1} - u^n}{\tau} v  d {\bm x} & +
  \int_{\Omega} k({\bm x}) \grad u^{n+1} \grad v d {\bm x} \\
  & +
  \int_{\partial \Omega} g({\bm x}) u^{n+1} v dx +
  p^{n+1} \int_{\Omega} u^{n} v d {\bm x} \\
  & = 
  \int_{\Omega} f({\bm x},t^{n+1}) v d {\bm x},
  \quad \forall v \in V,   
  \quad n = 0,1, ..., N-1 .
\end{split}
\end{equation}
In this case, the lower coefficient of the parabolic equation
is taken at the upper time level, whereas
the approximate solution $u({\bm x},t)$ is treated at the previous time level.
Let us consider the solution procedure of the problem  (\ref{eq:7}), (\ref{eq:9}), (\ref{eq:10}) in detail.

For the approximate solution at the new time level $u^{n+1}$, we introduce the following decomposition
\cite{samarskii2007numerical,borukhov2000numerical}:
\begin{equation}\label{eq:11}
  u^{n+1}({\bm x})  = y^{n+1}({\bm x})  + p^{n+1} w^{n+1}({\bm x}) .
\end{equation}
To find $y^{n+1}({\bm x})$, we employ the equation
\begin{equation}\label{eq:12}
\begin{split}
  \int_{\Omega} \frac{y^{n+1} - u^n}{\tau} v  d {\bm x} & +
  \int_{\Omega} k({\bm x}) \grad y^{n+1} \grad v d {\bm x} +
  \int_{\partial \Omega} g({\bm x}) y^{n+1} v dx \\
  & = 
  \int_{\Omega} f({\bm x},t^{n+1}) v d {\bm x},
  \quad \forall v \in V,   
  \quad n = 0,1, ..., N-1 .
\end{split}
\end{equation}
The function $w^{n+1}({\bm x})$ is determined from
\begin{equation}\label{eq:13}
\begin{split}
  \int_{\Omega} \frac{w^{n+1}}{\tau} v  d {\bm x} & +
  \int_{\Omega} k({\bm x}) \grad w^{n+1} \grad v d {\bm x} +
  \int_{\partial \Omega} g({\bm x}) w^{n+1} v dx \\
  & = - \int_{\Omega} u^{n} v d {\bm x} ,
  \quad \forall v \in V,   
  \quad n = 0,1, ..., N-1 .
\end{split}
\end{equation}
Using the decomposition (\ref{eq:11})--(\ref{eq:12}), equation (\ref{eq:10}) holds automatically for any $p^{n+1}$. 

To evaluate $p^{n+1}$, we apply the condition (\ref{eq:9}) (or (\ref{eq:8})). 
The substitution of (\ref{eq:11}) into (\ref{eq:9}) yields
\begin{equation}\label{eq:14}
  p^{n+1} =  \frac{1}{w^{n+1}({\bm x}^*)} (\varphi^{n+1} - y^{n+1}({\bm x}^*)) .
\end{equation}
The fundamental point of applicability of this algorithm is associated with the condition $w^{n+1}({\bm x}^*) \neq 0$.
The auxiliary function $w^{n+1}({\bm x})$ is determined from the grid elliptic equation (\ref{eq:13}).
The property of having fixed sign for $w^{n+1}({\bm x})$ is followed, in particular, from the same property of the solution at the previous time level $u^{n}({\bm x})$.
Such constraints on the solution can be provided by the corresponding restrictions on
the input data of the inverse problem. In any case, this problem requires special
and careful consideration. In this paper, we assume that the constraint $w^{n+1}({\bm x}^*) \neq 0$ is satisfied.

In solving problem (\ref{eq:8})--(\ref{eq:10}), instead of (\ref{eq:14}), we have
\begin{equation}\label{eq:15}
  p^{n+1} =  \frac{1}{\int_{\Omega} w^{n+1} \omega({\bm x})  d {\bm x}} \left (\varphi^{n+1} - \int_{\Omega} y^{n+1} \omega({\bm x})  d {\bm x} \right ) 
\end{equation} 
under the condition that
\[
\int_{\Omega} w^{n+1} \omega({\bm x})  d {\bm x} \neq 0 . 
\]
In this case, additional restrictions are formulated on the function $\omega({\bm x})$, e.g., its fixed sign in $\Omega$. 

Thus, the computational algorithm for solving the inverse problem
(\ref{eq:1})--(\ref{eq:4}) 
(or (\ref{eq:1})--(\ref{eq:3}), (\ref{eq:5})) based on
the linearized scheme (\ref{eq:7}), (\ref{eq:8}), (\ref{eq:10}) (or (\ref{eq:7}),(\ref{eq:9}), (\ref{eq:10}))
involves the solution of two standard grid
elliptic equations for the auxiliary functions $y^{n+1}({\bm x})$ (equation (\ref{eq:12})) 
and $w^{n+1}({\bm x})$ (equation (\ref{eq:13})), the further evaluation of
$p^{n+1}$ from (\ref{eq:15}) (or (\ref{eq:14})), and the final calculation $u^{n+1}({\bm x})$  from the relation (\ref{eq:11}).

\section{Scheme of the second-order accuracy} 
\label{sec:4}

The nonlinear inverse problem (\ref{eq:1})--(\ref{eq:4}) is characterized by a quadratic nonlinearity.
When using the scheme with linearization (\ref{eq:7}), (\ref{eq:8}), (\ref{eq:10}), the nonlinear term
is approximated with the first order with respect to $\tau$.
It is possible to apply the linearized scheme of second order. Let us consider the approximation
\[
  a(t^{n+1/2}) b(t^{n+1/2}) = \frac{1}{2} a(t^{n+1}) b(t^{n}) + \frac{1}{2} a(t^{n}) b(t^{n+1}) + O(\tau^2) . 
\]
Approximation of equation (\ref{eq:1}) with the boundary conditions (\ref{eq:2})
using the Crank-Nicolson scheme yields the linearized scheme
\begin{equation}\label{eq:16}
\begin{split}
  \int_{\Omega} \frac{u^{n+1} - u^n}{\tau} v  d {\bm x} & +
  \frac{1}{2} \int_{\Omega} k({\bm x}) \grad u^{n+1} \grad v d {\bm x} +
  \frac{1}{2} \int_{\partial \Omega} g({\bm x}) u^{n+1} v dx \\
  & + \frac{1}{2} \int_{\Omega} k({\bm x}) \grad u^{n} \grad v d {\bm x} +
  \frac{1}{2} \int_{\partial \Omega} g({\bm x}) u^{n} v dx \\
  & + \frac{1}{2} p^{n+1} \int_{\Omega} u^{n} v d {\bm x} +
  \frac{1}{2} p^{n} \int_{\Omega} u^{n+1} v d {\bm x} \\
  & = 
  \int_{\Omega} f({\bm x},t^{n+1}) v d {\bm x},
  \quad \forall v \in V,   
  \quad n = 0,1, ..., N-1 .
\end{split}
\end{equation} 
The scheme (\ref{eq:7}), (\ref{eq:8}), (\ref{eq:16}) belongs to the class of linearized schemes.
In comparison with the scheme (\ref{eq:7}), (\ref{eq:8}), (\ref{eq:10}),
it has a higher order of accuracy in time.

To implement (\ref{eq:7}), (\ref{eq:8}), (\ref{eq:16}), we again use the decomposition (\ref{eq:11}).
In this case, for $y^{n+1}({\bm x})$, we have 
\begin{equation}\label{eq:17}
\begin{split}
  \int_{\Omega} \frac{y^{n+1} - u^n}{\tau} v  d {\bm x} & +
  \frac{1}{2}\int_{\Omega} k({\bm x}) \grad y^{n+1} \grad v d {\bm x}  +
  \frac{1}{2}\int_{\partial \Omega} g({\bm x}) y^{n+1} v dx \\
  & + \frac{1}{2} \int_{\Omega} k({\bm x}) \grad u^{n} \grad v d {\bm x} \\
  & +
  \frac{1}{2} \int_{\partial \Omega} g({\bm x}) u^{n} v dx +
  \frac{1}{2} p^{n} \int_{\Omega} y^{n+1} v d {\bm x} \\
  & = 
  \int_{\Omega} f({\bm x},t^{n+1}) v d {\bm x},
  \quad \forall v \in V,   
  \quad n = 0,1, ..., N-1 .
\end{split}
\end{equation}
The auxiliary function $w^{n+1}({\bm x})$  is defined as the solution of the equation
\begin{equation}\label{eq:18}
\begin{split}
  \int_{\Omega} \frac{w^{n+1}}{\tau} v  d {\bm x} & +
  \frac{1}{2} \int_{\Omega} k({\bm x}) \grad w^{n+1} \grad v d {\bm x} \\
  & +
  \frac{1}{2} \int_{\partial \Omega} g({\bm x}) w^{n+1} v dx +
  \frac{1}{2} p^{n} \int_{\Omega} w^{n+1} v d {\bm x} \\
  & = - \frac{1}{2} \int_{\Omega} u^{n} v d {\bm x} ,
  \quad \forall v \in V,   
  \quad n = 0,1, ..., N-1 .
\end{split}
\end{equation}
Further, as in the case of the first-order scheme, we employ (\ref{eq:15}).

The Crank-Nicholson scheme for numerical solving the direct problems for
parabolic equations is not very often used in computational practice.
It is inferior to the fully implicit scheme in sense of conservation of monotonicity (fulfilment of
the maximum principle for the grid problem), it has poor asymptotic properties for
solving problems with large integration time, and it is not unconditionally SM-stable scheme
\cite{samarskiaei1995computational,vabishchevich2012sm}.
For this reason, it is appropriate to consider another variant of linearization of this inverse problem, 
where the  second-order approximation is applied only for the nonlinear term.
In this case, instead of (\ref{eq:10}), (\ref{eq:16}), we put
\begin{equation}\label{eq:19}
\begin{split}
  \int_{\Omega} \frac{u^{n+1} - u^n}{\tau} v  d {\bm x} & +
  \int_{\Omega} k({\bm x}) \grad u^{n+1} \grad v d {\bm x} +
  \int_{\partial \Omega} g({\bm x}) u^{n+1} v dx \\
  & + \frac{1}{2} p^{n+1} \int_{\Omega} u^{n} v d {\bm x} +
  \frac{1}{2} p^{n} \int_{\Omega} u^{n+1} v d {\bm x} \\
  & = 
  \int_{\Omega} f({\bm x},t^{n+1}) v d {\bm x},
  \quad \forall v \in V,   
  \quad n = 0,1, ..., N-1 .
\end{split}
\end{equation} 
The numerical implementation of the scheme (\ref{eq:7}), (\ref{eq:8}), (\ref{eq:19})
is performed in the standard way using the decomposition (\ref{eq:11}).

\section{Numerical examples} 
\label{sec:5}

\begin{figure}[htp]
  \begin{center}
    \includegraphics[scale = 0.35] {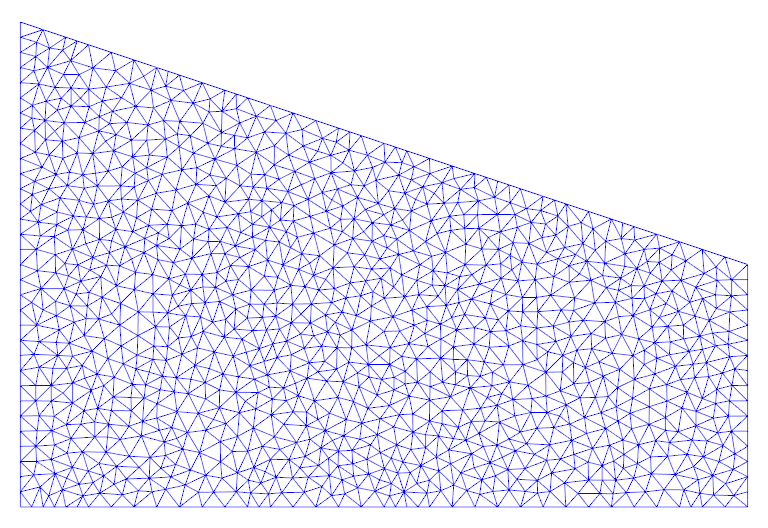}
	\caption{The computational grid}
	\label{fig:1}
  \end{center}
\end{figure}

\begin{figure}[htp]
  \begin{center}
    \includegraphics[scale = 0.4] {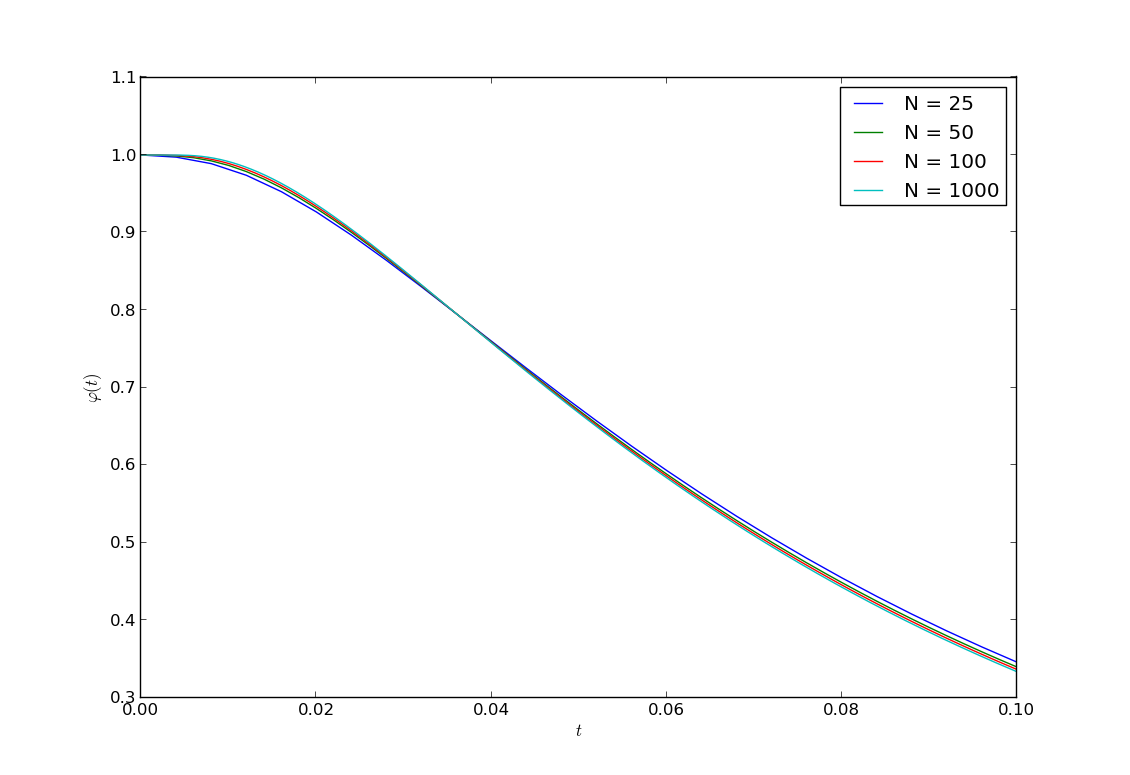}
	\caption{The solution of the direct problem at the point of observation}
	\label{fig:2}
  \end{center}
\end{figure}

\begin{figure}[htp]
  \begin{center}
    \includegraphics[scale = 0.275] {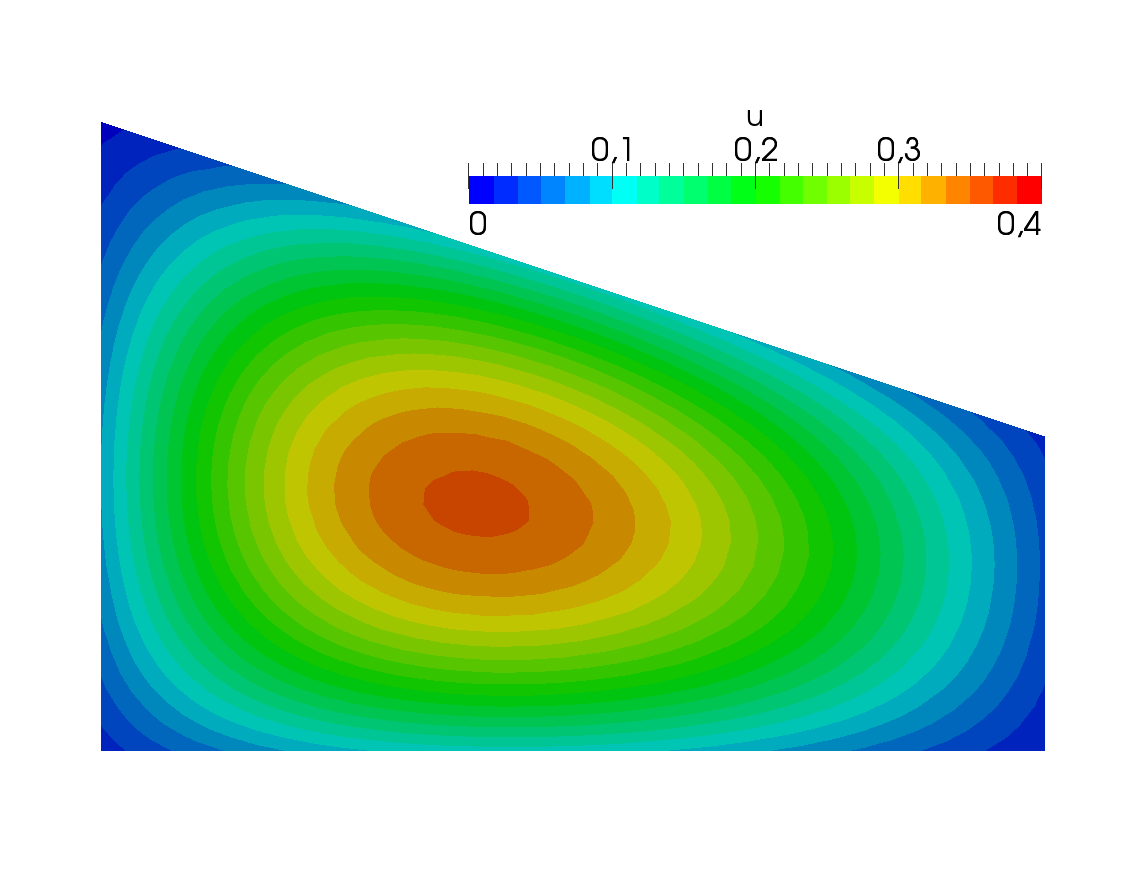}
	\caption{The solution of the direct problem at $t=T$}
	\label{fig:3}
  \end{center}
\end{figure}

To demonstrate possibilities of the above linearization schemes for solving
the problem of the identification of the lower coefficient of the parabolic equation,
we consider a 2D model problem. In the examples below, we put
\[
  k({\bm x}) = 1, 
  \quad f({\bm x},t) = 0, 
  \quad u_0({\bm x}) = 1,   
  \quad  {\bm x} \in \Omega,
    \quad g({\bm x}) = 10,  
  \quad  {\bm x} \in \partial \Omega.
\]
The problem is considered on a triangular grid, which consists of 1,180 nodes (2,230 triangles) and
is shown in Fig.\ref{fig:1}. Here $\Omega$ is the trapezoid with the vertices coordinates $(0,0), (0,1), (1.5, 0.5), 1.5, 0)$.
The calculations were carried out for $T = 0.1$. The coefficient $p(t)$ is taken in the form 
\begin{equation}\label{eq:20}
  p(t) = \left \{
  \begin{array}{cc}
  1000 t, & \quad 0 < t \leq \frac{1}{2} T, \\
  0, & \quad \frac{1}{2} T < t \leq T . \\   
  \end{array}  
  \right . 
\end{equation} 
The solution of the direct problem (\ref{eq:1})--(\ref{eq:3}) at the observation point 
is depicted in Fig.\ref{fig:2}. It was otained using the fully implicit scheme with different time steps.
The solution at the final time moment is presented in  Fig.\ref{fig:3}. 


\begin{figure}[htp]
  \begin{center}
    \includegraphics[scale = 0.4] {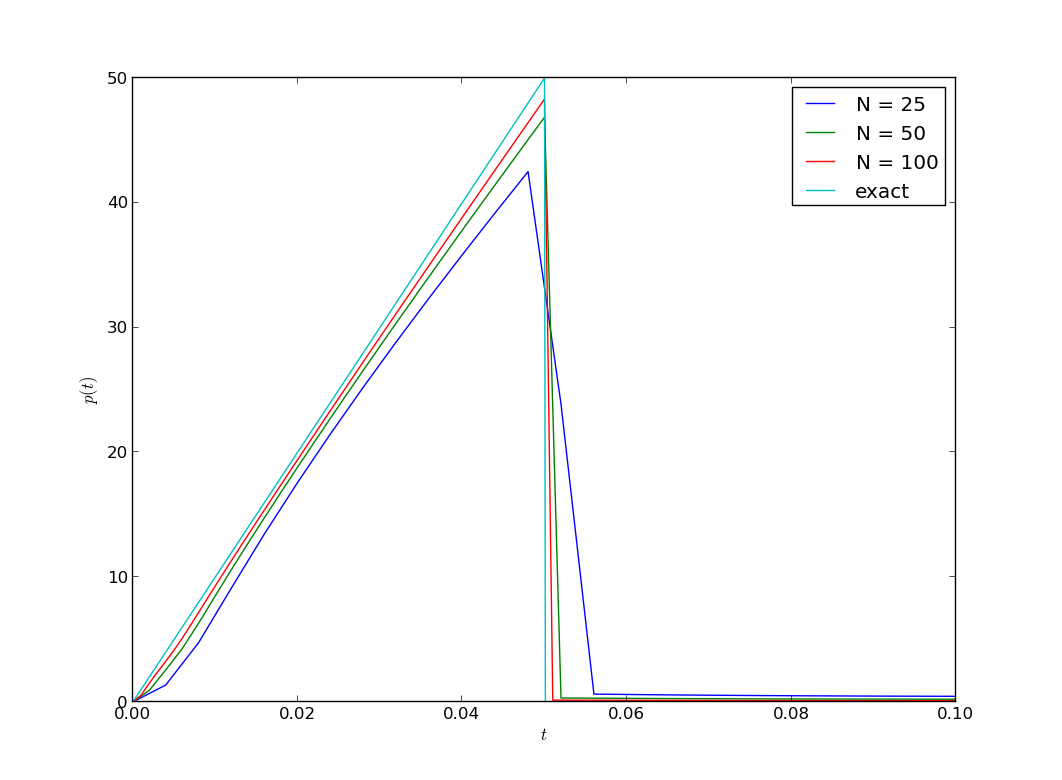}
	\caption{The solution of the inverse problem}
	\label{fig:4}
  \end{center}
\end{figure}

The results of solving the inverse problem with variuos grids in time
are shown in Fig.\ref{fig:4}. The solution of the direct problem obtained with $N = 1000$
is used as the input data (the function $\varphi(t)$ in the  condition (\ref{eq:5})).
It is easy to see that the approximate solution of the inverse problem
converges  with decreasing the time step. These results were obtained using the first-order
scheme (\ref{eq:7}), (\ref{eq:9}), (\ref{eq:10}). 

\begin{figure}[htp]
  \begin{center}
    \includegraphics[scale = 0.4] {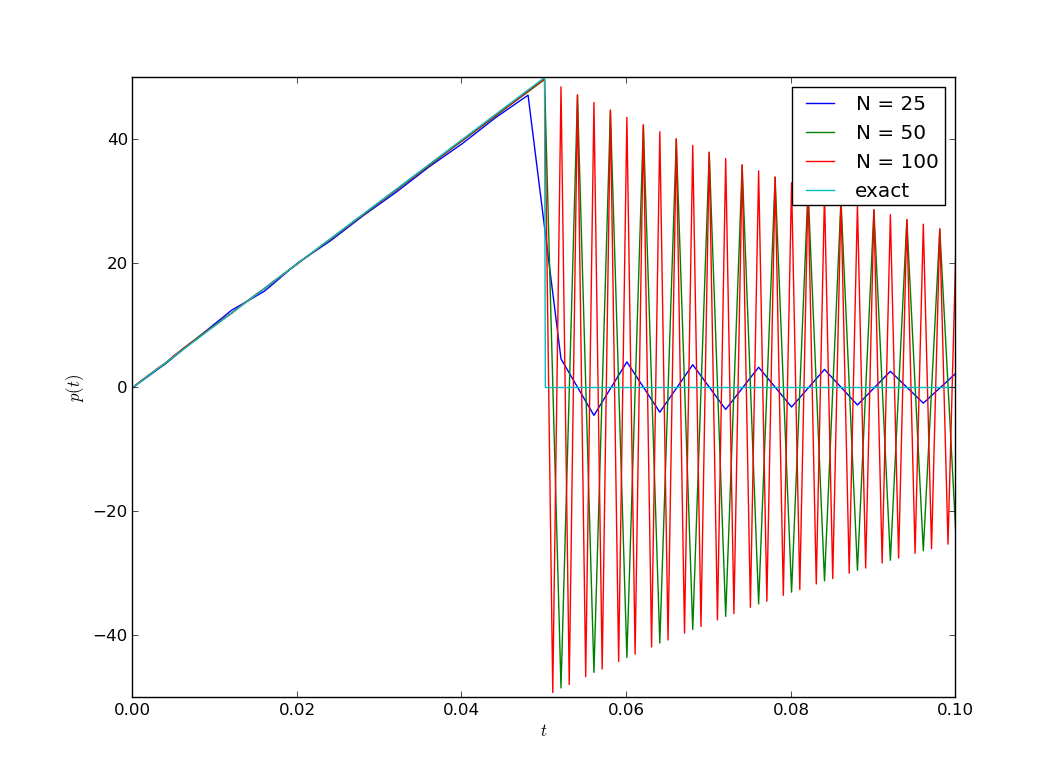}
	\caption{The solution of the inverse problem (the second-order scheme)}
	\label{fig:5}
  \end{center}
\end{figure}

Numerical results obtained for the above problem
using the second-order scheme (\ref{eq:7}), (\ref{eq:9}), (\ref{eq:16})
are shown in Fig.\ref{fig:5}. For the discontinuous right-hand side (\ref{eq:20}), 
we observe characteristic wiggles of the identified coefficient.
Such oscillations of the approximate solution are typical for the scheme
(\ref{eq:7}), (\ref{eq:9}), (\ref{eq:19}).

\begin{figure}[htp]
  \begin{center}
    \includegraphics[scale = 0.4] {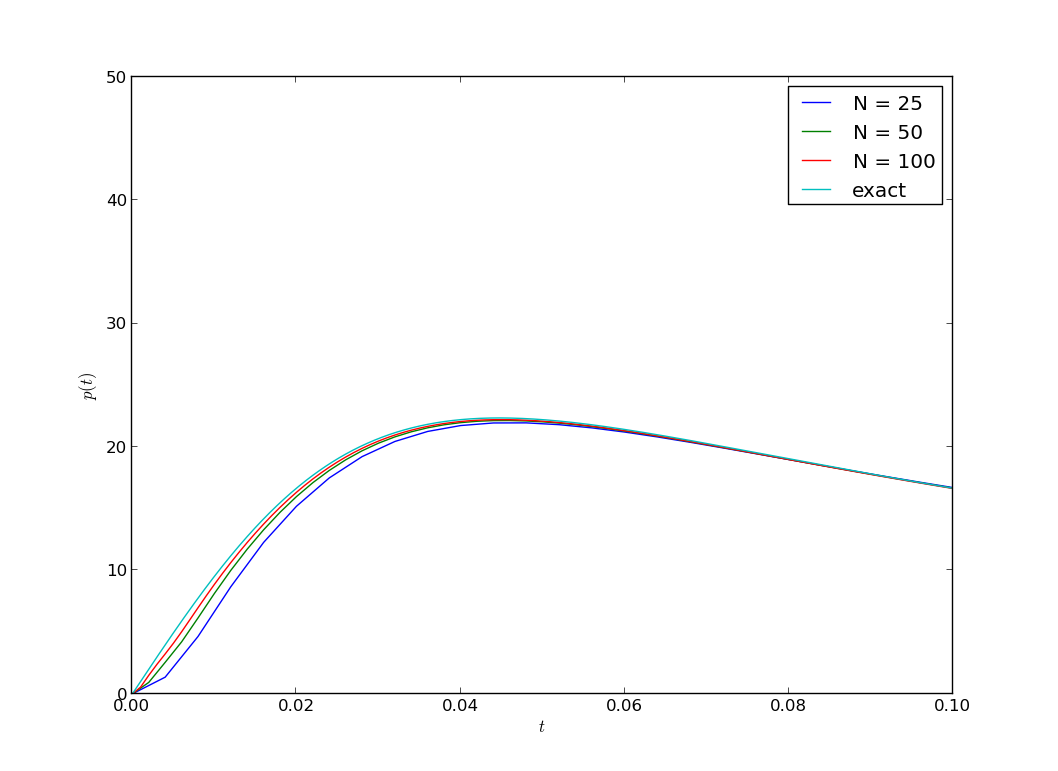}
	\caption{The scheme of first order (smooth solution)}
	\label{fig:6}
  \end{center}
\end{figure}

\begin{figure}[htp]
  \begin{center}
    \includegraphics[scale = 0.4] {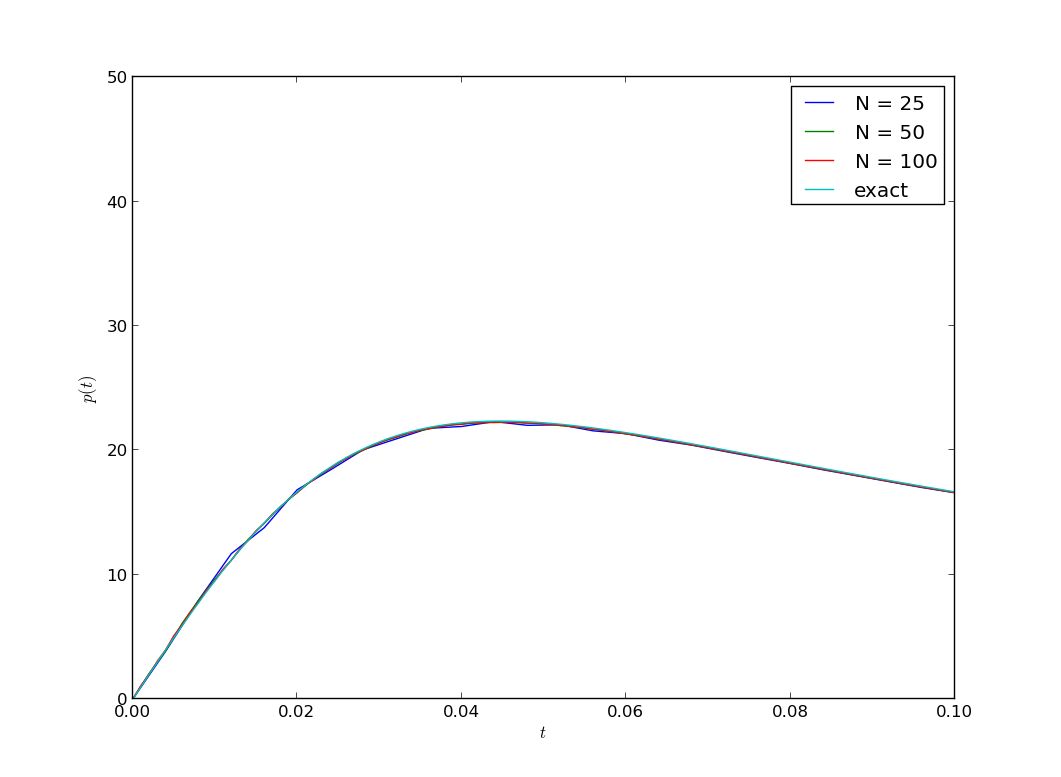}
	\caption{The scheme of second order (smooth solution)}
	\label{fig:7}
  \end{center}
\end{figure}

If the desired solution (the coefficient $p(t)$) is smooth, then the effect of using
the second-order approximation is clearly expressed. As an example, we present the results of numerical solving
the inverse problem, where the lower coefficient (the exact solution) has the form
\[
 p(t) = \frac{1000 t}{1+ 500 t^2} .
\] 
The approximate solution obtained via the first-order scheme (\ref{eq:7}), (\ref{eq:9}), (\ref{eq:10})  
is shown in Fig.\ref{fig:6}, whereas Fig.\ref{fig:7} demonstrates the computations conducted by means of 
the  second-order scheme (\ref{eq:7}), (\ref{eq:9}), (\ref{eq:16}).



\begin{thebibliography}{10}
\expandafter\ifx\csname url\endcsname\relax
  \def\url#1{\texttt{#1}}\fi
\expandafter\ifx\csname urlprefix\endcsname\relax\def\urlprefix{URL }\fi
\expandafter\ifx\csname href\endcsname\relax
  \def\href#1#2{#2} \def\path#1{#1}\fi

\bibitem{alifanov2011inverse}
O.~M. Alifanov, Inverse Heat Transfer Problems, Springer, 2011.

\bibitem{tarantola1987inverse}
A.~Tarantola, Inverse problem theory: methods for data fitting and model
  parameter estimation, Elsevier, 1987.

\bibitem{aster2011parameter}
R.~C. Aster, B.~Borchers, C.~H. Thurber, Parameter Estimation and Inverse
  Problems, Elsevier Science, 2011.

\bibitem{kirsch1996introduction}
A.~Kirsch, An introduction to the mathematical theory of inverse problems,
  Springer, 1996.

\bibitem{tikhonov1977solutions}
A.~N. Tikhonov, V.~Y. Arsenin, Solutions of ill-posed problems, Winston, 1977.

\bibitem{engl2000regularization}
H.~W. Engl, M.~Hanke, A.~Neubauer, Regularization of inverse problems, Kluwer
  Academic Publishers, 2000.

\bibitem{lavrentev1986ill}
M.~M. Lavrent'ev, V.~G. Romanov, S.~P. Shishatskii, Ill-posed Problems of
  Mathematical Physics and Analysis, American Mathematical Society, 1986.

\bibitem{isakov1998inverse}
V.~Isakov, Inverse problems for partial differential equations, Springer, 1998.

\bibitem{express2000methods}
A.~I. Prilepko, D.~G. Orlovsky, I.~A. Vasin, Methods for Solving Inverse
  Problems in Mathematical Physics, Marcel Dekker, Inc, 2000.

\bibitem{belov2002inverse}
Y.~Y. Belov, Inverse problems for partial differential equations, VSP, 2002.

\bibitem{maksimov2002dynamical}
V.~I. Maksimov, Dynamical inverse problems of distributed systems, VSP, 2002.

\bibitem{vogel2002computational}
C.~R. Vogel, Computational Methods for Inverse Problems, no.~10, Society for
  Industrial and Applied Mathematics, 2002.

\bibitem{samarskii2007numerical}
A.~A. Samarskii, P.~N. Vabishchevich, Numerical Methods for Solving Inverse
  Problems of Mathematical Physics, De Gruyter, 2007.

\bibitem{prilepko1985determination}
A.~I. Prilepko, D.~G. Orlovskii, Determination of evolution parameter of an
  equation, and inverse problems in mathematical physics, part {I} and {II},
  Differ. Equat. 21~(1, 4) (1985) 119--129, 694--701, in Russian.

\bibitem{prilepko1987solvability}
A.~I. Prilepko, V.~V. Soloev, Solvability of the inverse boundary value problem
  of finding a coefficient of a lower order term in a parabolic equation,
  Differ. Equat. 23~(1) (1987) 136--143, in Russian.

\bibitem{macbain1986existence}
J.~A. MacBain, J.~B. Bednar, Existence and uniqueness properties for the
  one-dimensional magnetotellurics inversion problem, Journal of mathematical
  physics 27 (1986) 645--649.

\bibitem{cannon1990inverse}
J.~R. Cannon, Y.~Lin, An inverse problem of finding a parameter in a
  semi-linear heat equation, Journal of Mathematical Analysis and Applications
  145~(2) (1990) 470--484.

\bibitem{wang1989finite}
S.~Wang, Y.~Lin, A finite-difference solution to an inverse problem for
  determining a control function in a parabolic partial differential equation,
  Inverse problems 5~(4) (1989) 631--640.

\bibitem{cannon1994numerical}
J.~R. Cannon, Y.~Lin, S.~Xu, Numerical procedures for the determination of an
  unknown coefficient in semi-linear parabolic differential equations, Inverse
  Problems 10~(2) (1994) 227--243.

\bibitem{dehghan2005identification}
M.~Dehghan, Identification of a time-dependent coefficient in a partial
  differential equation subject to an extra measurement, Numerical Methods for
  Partial Differential Equations 21~(3) (2005) 611--622.

\bibitem{ye2007stability}
C.~Ye, Z.~Sun, On the stability and convergence of a difference scheme for an
  one-dimensional parabolic inverse problem, Applied mathematics and
  computation 188~(1) (2007) 214--225.

\bibitem{wang2012inverse}
W.~Wang, B.~Han, M.~Yamamoto, Inverse heat problem of determining
  time-dependent source parameter in reproducing kernel space, Nonlinear
  Analysis: Real World Applications 14~(1) (2013) 875--887.

\bibitem{dehghan2001numerical}
M.~Dehghan, Numerical methods for two-dimensional parabolic inverse problem
  with energy overspecification, International journal of computer mathematics
  77~(3) (2001) 441--455.

\bibitem{daoud2005splitting}
D.~S. Daoud, D.~Subasi, A splitting up algorithm for the determination of the
  control parameter in multi dimensional parabolic problem, Applied mathematics
  and computation 166~(3) (2005) 584--595.

\bibitem{quarteroni2008numerical}
A.~Quarteroni, A.~Valli, Numerical approximation of partial differential
  equations, Springer, 2008.

\bibitem{brenner2008mathematical}
S.~C. Brenner, L.~R. Scott, The mathematical theory of finite element methods,
  Springer, 2008.

\bibitem{borukhov2000numerical}
V.~T. Borukhov, P.~N. Vabishchevich, Numerical solution of the inverse problem
  of reconstructing a distributed right-hand side of a parabolic equation,
  Computer physics communications 126~(1) (2000) 32--36.

\bibitem{samarskiaei1995computational}
A.~A. Samarskii, P.~N. Vabishchevich, Computational heat transfer, Wiley, 1995.

\bibitem{vabishchevich2012sm}
P.~N. Vabishchevich, {SM}-stability of operator-difference schemes,
  Computational Mathematics and Mathematical Physics 52~(6) (2012) 887--894.

\end{thebibliography}
\end{document}